\documentclass[5p,doublecolumn]{elsarticle}
\usepackage{atbegshi}
\AtBeginDocument{\AtBeginShipoutNext{\AtBeginShipoutDiscard}}                                   
\usepackage{mathrsfs}
\usepackage{tikz}
\usetikzlibrary{decorations.pathmorphing}
\usepackage{amssymb}
\usepackage{amsmath}               
\usepackage{amsfonts}              
\usepackage{graphicx}          
\usepackage{subfigure}
\newcounter{ctr}
\newcounter{ctr2}
\newcounter{ctr3}
\newcounter{ctr4}
\newcounter{ctr5}
\newtheorem{coro}[ctr]{Corollary}
\newtheorem{them}[ctr2]{Theorem}
\newtheorem{defi}[ctr3]{Definition}
\newtheorem{lemm}[ctr4]{Lemma}
\newtheorem{prop}[ctr5]{Proposition}

\newenvironment{pf}{\paragraph{Proof:}}{\hfill$\square$}

\begin{document}
\begin{frontmatter}

\title{On ~Robust ~Stability ~of ~Switched ~Systems~in~the~Context~of~Filippov~Solutions} 


\author[First]{Mohamadreza Ahmadi}\ead{mrahmadi@utexas.edu},    
\author[Second]{Hamed Mojallali}\ead{mojallali@guilan.ac.ir},               
\author[Third]{Rafael Wisniewski}\ead{raf@es.aau.dk}  

\address[First]{Institute for Computational Engineering and Sciences (ICES), University of Texas at Austin, Peter O'Donnell, Jr. Building, 201 E 24th St, Austin, TX 78712.}
\address[Second]{Department of Electrical Engineering, Faculty of Engineering, University of Guilan, Rasht, Iran,  P.O. Box: 41635-3756, Postal Code: 41996-13769.}
\address[Third]{Department of Electronic Systems, Automation and Control, Aalborg University, Fredrik Bajers Vej 7 C, 9220, Aalborg East, Denmark.}


\begin{abstract}                          
The stability problem of a class of nonlinear switched systems defined on compact sets with state-dependent switching   is considered. Instead of the Carath\'eodory solutions, the general Filippov solutions are studied. This encapsulates solutions with infinite switching in finite time and sliding modes in the neighborhood of the switching surfaces.  In this regard, a Lyapunov-like stability theorem, based on the theory of differential inclusions, is formulated. Additionally, the results are extended to switched systems with simplical uncertainty. It is also demonstrated that, for the special case of polynomial switched systems defined on semi-algebraic sets, stability analysis can be checked based on sum of squares programming techniques.
\end{abstract}

\end{frontmatter}

\pagenumbering{arabic}

\section{Introduction}
A plethora of systems encountered in engineering and nature give rise to mathematical models which encompass both discrete and continuous dynamics. Conventionally, these systems are identified by a family of indexed differential or difference equations describing each subsystem and a switching rule between them. This rich family of systems is referred to switched or more generally hybrid systems. Due to their ubiquitous nature, a significant amount of literature has emerged on investigating their real world applications; e.g.,~\cite{weld}. Furthermore, the analysis of switched and hybrid systems has received tremendous attention \cite{forte2016robust,papusha2016automata,murti2013sum,kundu2016generalized,Ali201727,7206541,WANG201694}.

However despite their prevalence, the stability issue of switched systems has not yet been completely resolved \cite{ant09}, \cite{siam07}. Several interesting phenomena arise when dealing with such systems; to name but a few, even if all the subsystems are exponentially stable, one cannot guarantee the stability of the overall system \cite{bran}. Conversely, an appropriate switching law may contribute to stability even when all subsystems are unstable \cite{liberzon}.  Besides, a switched system can exhibit chaotic dynamics which further exacerbates stability problems \cite{chaos93}. Still, the notion of stability for switched systems is contingent on the  type of solutions considered~\cite{Georgescu20121985}. This can be exemplified as a switched system with stable Carath\'eodory solutions, may possess divergent Filippov solutions (see Example 5 in \cite{lewis}). Therefore, stability of Carath\'eodory solutions does not imply the overall stability of the corresponding switched system.

It has been demonstrated that exploiting the theory of differential inclusions is a promising methodology for describing and analyzing the dynamics of switched and hybrid systems.  In a seminal contribution, Botchkaref and Tripakis \cite{verifi} studied the verification of a class of hybrid systems characterized by linear differential inclusions. Aubin \textit{et. al.} \cite{impulse} brought forward conditions to determine viable or invariant states of a hybrid system described by impulsive differential inclusions.  Margaliot and Liberzon \cite{margaliot} proposed Lie-algebraic stability conditions for a relaxed differential inclusion representing a switched system.  Goebel \textit{et. al.} \cite{goebel} formulated asymptotic stability conditions for hybrid dynamical systems defined by differential and difference inclusions.  Leth and Wisniewski \cite{lewis} applied the theory of differential inclusions and suggested Lyapunov-like stability theorem for switched systems defined on polyhedral sets. Following the same trend, Ahmadi \textit{et. al.} \cite{ahmadi2012}, \cite{ahmadi2014guaranteed} presented a robust $H_\infty$ controller synthesis scheme for the latter class of systems subject to uncertainty.

This study is predominantly motivated by \cite{pranja} and \cite{lewis}. In \cite{pranja}, Prajna and Papachristodoulou efficaciously put forth stability analysis tools for a class of hybrid systems using sum of squares (SOS) techniques; but, the authors did not advance a corresponding well-founded stability theorem. Moreover, the solutions implicitly considered in \cite{pranja} are in the sense of Carath\'eodory, which connote the exclusion of solutions with infinite switching in finite time from the analysis. On the other hand, \cite{lewis} is concerned with stability of switched systems defined on polyhedral sets (piecewise affine systems) in the framework of Filippov solutions (see also the intriguing discussions maintained in \cite{Pogromsky20032139} and \cite{Heemels20083079}); however, it does not provide any computational tools for determining stability. In the present paper, the stability results delineated in \cite{lewis} are generalized to nonlinear switched systems defined on compact sets. This generalization is carried out established upon the theoretical results from  differential inclusions. Furthermore,  the robust stability problem of switched systems with simplical uncertainty is addressed. Subsequently, in order to provide the means of computationally efficient analysis, we propose sufficient conditions based on SOS programming for the suggested stability theorems. Simulation results are also supplemented which corroborate the theoretical analyses given is the paper.

The framework of this paper proceeds as follows. The notations and some preliminary mathematical discussions adopted in this study are limned in Section \ref{secpre}. The main contributions of this paper are outlined in Section~\ref{secmain2}. The proposed methodologies are elucidated in Section~\ref{secexam} via a simulation example. Finally, Section~\ref{seccon} concludes the paper.

\section{Notations and Preliminaries} \label{secpre}

The notations employed in this paper are relatively straightforward. The set of non-negative real numbers $[0,\infty)$ is denoted by $\mathbb{R}_{\ge 0}$. The Euclidean vector norm on $\mathbb{R}^n$ is designated by $\| \cdot \|$, the inner product by $\langle \cdot \rangle$, and the closed ball of radius $\epsilon$ in $\mathbb{R}^n$ centered at origin by $\mathcal{B}_{\epsilon}^n$. Let $\mathcal{P}(x)$ account for the ring of polynomial functions $p: \mathbb{R}^n \to \mathbb{R}$ over $\mathbb{R}$ in the variable $x$ and $\mathcal{P}_{sos}(x)\subset\mathcal{P}(x)$ the subset of polynomials with an SOS decomposition; i.e, $p(x) \in \mathcal P_{sos}(x)$ if and only if there are $p_i(x) \in \mathcal P(x), i \in \{1, \ldots ,k\}$ such that $p(x) = p_1^2(x) + \cdots +p_k^2(x)$. We denote the interior of a compact set $K$ by $int(K)$, and the boundary of $K$ by $bd(K)$; then, $K=int(K) \cup bd(K)$. The closed convex hull of the set $K$ is denoted by $co(K)$, and the set of all subsets of $K$ (power set of $K$) is represented by $2^K$.

\subsection{Switched Systems Defined on Compact Sets}

In this study, we partition the state space by a family of closed sets. On this family, we impose certain regularity conditions, as delineated in the definition below.

\begin{defi}[Nice Covering]
Let $\mathcal{X}$ be a compact subset set of an Euclidean space $E$ with
\begin{equation} \label{cscs}
cl(int(\mathcal{X} )) = \mathcal{X}. 
\end{equation}
 A family $\mathscr{X} = \{X_i\}_{i\in I}$ of subsets of $\mathcal{X}$ is a nice covering of
$\mathcal{X}$ if and only if
 \begin{enumerate}
\item $\mathscr{X}$ is a covering of $\mathcal{X}$, i.e., $\mathcal{X} =  \bigcup_{i\in I} X_i$, and \\$cl(int(X_i)) = X_i$ for all $i \in I$;
\item $X_i \cap X_j =bd(X_i) \cup bd(X_j)$ for all $(i,j)\in I^2$;
\item $\mathscr{X}$ is locally finite, i.e., each point of $\mathcal{X}$ has an open neighborhood intersecting only finitely many elements of $\mathscr{X}$;
\item for any $x\in bd(X_i)$ and $v\in E$, there is $a>0$ and $j\in I$ such
that $x+vt \in X_j$ for all $t \in [0,a)$.
\end{enumerate}
\end{defi}

In this study, we consider a class of $n$-dimensional nonlinear switched systems $\mathcal{S}= \lbrace \mathcal{X},\mathscr{X},I,\mathscr{F} \rbrace$, wherein $\mathcal{X}$ is a compact set representing the state-space that satisfies~\eqref{cscs}, $\mathscr{X}= \lbrace X_i \rbrace_{i \in I}$ is a nice covering of $\mathcal{X}$ with index set $I$, and $\mathscr{F}=\lbrace F_i \rbrace_{i \in I}$ a family of smooth vector fields. Each function $F_i$ is defined on an open neighborhood of origin $U_i$  ($F_i : U_i \to \mathbb{R}^n$). 

For a nice covering $\mathscr{X} = \{ X_i \}_{i \in I}$ with index set $i \in I$, we define $\tilde{I}= \left \{ (i,j) \in I^2 \mid X_i \cap X_j \neq \emptyset,~i \neq j   \right\}$, the set of index pairs which determines the partitions with non-empty intersections. Remark that partitioning by polyhedral sets assures that this latter property is satisfied. This is the case when considering switched systems defined on polyhedral sets; e.g., piecewise affine systems. We shall say that a switching has occurred whenever a trajectory passes some boundary (switching surface).

 The global dynamics is described by the following differential inclusions
\begin{equation} \label{eq1}
\dot{x}(t) \in \mathcal{F} \left(x(t)\right)
\end{equation}
\begin{equation} \label{eq2}
\dot{x}(t) \in \mathcal{F}^c \left(x(t)\right),
\end{equation}
where the set-valued maps $\mathcal{F}$ and $\mathcal{F}^c$ are defined by
\begin{equation} \label{eq3}
\mathcal{F} : \mathcal{X} \to 2^\mathcal{X}; x \mapsto  \left \lbrace v \in \mathbb{R}^n \mid v=F_i(x)  \text{  if  }  x \in X_i \right \rbrace
\end{equation}
\begin{equation} \label{eq4}
\mathcal{F}^c : \mathcal{X} \to 2^\mathcal{X}; x \mapsto  co \left(\mathcal{F}(x)\right),
\end{equation}
 The choice of whether the dynamics is modeled by \eqref{eq1} or \eqref{eq2} depends on the nature of the motion to be considered (see Fig.~\ref{fig1}). Pertaining to the solutions of discontinuous and switched dynamical systems, the interested reader is referred to the didactic review in \cite{cort}.

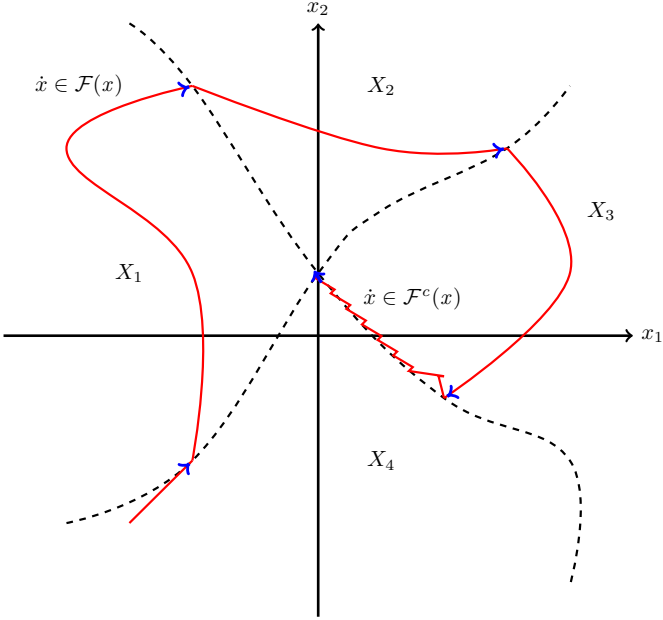
\begin{figure}[t] \label{fig1}
\begin{center}
\resizebox{9cm}{!}{
\begin{tikzpicture}

  \begin{scope} [line width=1.2]
    \draw[->] (-5,0) -- (5,0) node[right] {$x_1$};
    \draw[->] (0,-4.5) -- (0,5) node[above] {$x_2$};

  \end{scope}

\draw[dashed,line width=1]  plot[smooth, tension=.7] coordinates {(-4,-3)   (-2,-2) (0,1) (1,2)  (3,3)(4,4)};
\draw[dashed,line width=1]  plot[smooth, tension=.7] coordinates {(-3,5)(-2,4)  (0,1) (2,-1) (4,-2) (4,-4)};
[line width=1.2 ]
\draw[red, line width=1]  plot[smooth, tension=.7] coordinates {(-3,-3)  (-2,-2)};
\draw[red, line width=1]  plot[smooth, tension=.7] coordinates {(-2,-2)  (-2,1)   (-4,3) (-2,4) };
\draw[red, line width=1]  plot[smooth, tension=.7] coordinates {(-2,4) (1,3)(3,3)};
\draw[red, line width=1]  plot[smooth, tension=.7] coordinates {(3,3)(4,1) (2,-1)};

\draw[red,line width=1] [decorate, decoration=saw] plot[smooth, tension=.7] coordinates {(2,-1) (0,1)};
\draw [-to,shorten >=2pt,blue,ultra thick] (-2.05,-2.05) -- (-2,-2);
\draw [-to,shorten >=2pt,blue,ultra thick] (-2.05,3.985) -- (-2,4);
\draw [-to,shorten >=2pt,blue,ultra thick] (2.95,2.99) -- (3,3);
\draw [-to,shorten >=2pt,blue,ultra thick] (2.05,-0.97) -- (2,-1);
\draw [-to,shorten >=2pt,blue,ultra thick] (0.05,0.85) -- (-0.1,1.1);
\draw (-3,1) node {$X_1$};
\draw (1,4) node {$X_2$};
\draw (4.5,2) node {$X_3$};
\draw (1,-2) node {$X_4$};
\draw (-3.8,4) node {$\dot x \in\mathcal F( \mathit x )$};
\draw (1.5,0.6) node {$\dot x \in\mathcal F^\mathit c( \mathit x )$};
\end{tikzpicture}
}
\caption{The trajectories of a switched system. Notice that the motion follows either differential inclusions \eqref{eq1} or \eqref{eq2}.}  
\label{fig1}                                 
\end{center}                                 
\end{figure}
In the sequel, we apply the following notions from the theory of differential inclusions. For $Q \subset \mathbb{R}^n$, $T_{Q}(x)$ denotes the Bouligand's contingent cone\footnote{The Bouligand contingent cone $T_Q(x)$ to $Q$ at $x \in \mathcal{X}$ is the set of directions $v \in \mathcal{X}$ such that there exist sequences $h_n \to 0^{+}$ ($h_n \ge 0$) and $v_n \to v$ such that $x+h_n v_n \in \mathcal{X}$ for all $n \in \mathbb{N}$.} of $Q$ at $x \in Q$. If $Q$ is convex then $T_{Q}(x)$ is closure of the cone spanned by $Q-\{x\}$. In addition, if $x \in int(Q)$, we have $T_{Q}(x)=\mathbb{R}^n$ \cite{aubin}. The upper contingent derivative of a function $U : \mathbb{R}^n \to \mathbb{R}$ at $x_0$ in the direction $v_0$ is defined as
\begin{equation} \label{eq5}
D_{+}U(x_0)(v_0)=\lim_{\substack{h \to 0^{+} \\ v \to v_0}} \inf \left( \frac{U(x+hv)-U(x)}{h} \right)
\end{equation}
if $U$ is locally Lipschitzean then
\begin{equation}
D_{+}U(x_0)(v_0)=\lim_{{h \to 0^{+}}} \inf \left( \frac{U(x+hv_0)-U(x)}{h} \right)
\end{equation}
Note in particular that, if $U(x)$ is continuously differentiable it holds that
\begin{equation} \label{eq6}
D_{+}U(x_0)(v_0)=  \left \langle \frac{\partial{U}}{\partial{x}}, v \right \rangle
\end{equation}

\begin{prop} \label{p1}
The set-valued map $\mathcal{F}$ defined by \eqref{eq3} is upper semi-continuous. 
Correspondingly, $\mathcal{F}^c$ defined by \eqref{eq4} is an upper semi-continuous set valued map with non-empty, convex and compact values.
\end{prop}
\begin{pf}
For all $x \in int (X_i)$, $i \in I$, $\mathcal F(x) = F_i(x)$ is a one point set and since each $F_i$ is continuous, $\mathcal F$ is upper semi-continuous at any $x \in int (X_i)$. Furthermore, for any $x \in bd (X_i) \cup bd (X_{i+1}) \cup \cdots \cup bd (X_j)$, $\mathcal F(x) = \{ F_i(x),F_{i+1}(x),\ldots,F_j(x) \}$ is a multi-valued set. Because each $F_i$ is continuous, from the Weierstrass definition of continuity, it follows that for all $x \in X_i$ and $\epsilon > 0$ there exists a $\sigma_i>0$ such that $\forall \acute{x} \in x + \mathcal{B}_{\sigma_i}^n, \, F_i(\acute{x}) \in F_i(x) + \mathcal{B}_{\epsilon}^n$. To demonstrate that $\mathcal F$ is upper semi-continuous, it suffices to choose $\sigma = \min_{i} \sigma_i$. Then, for any $x \in \mathcal{X}$ and any $\epsilon > 0$ there exists a $\sigma (\epsilon,x) \le \epsilon$ such that $\forall \acute{x} \in x + \mathcal B_\sigma^n, \quad \mathcal F(\acute{x}) \subset \mathcal F(x) + \mathcal{B}_{\epsilon}^n$.   Additionally, because each of the maps $F_i, \, i \in I$, is continuous and $\mathcal F(x)$ is finite for all $x \in \mathcal{X}$ (finiteness of $\mathcal F$ follows from construction covering with local finiteness property); then, from Lemma 16 in p. 66, \cite{filip}, it follows that $\mathcal{F}^c$ is also upper semi-continuous.
\end{pf}

\begin{figure}
\centering
  \includegraphics[width=0.4\textwidth]{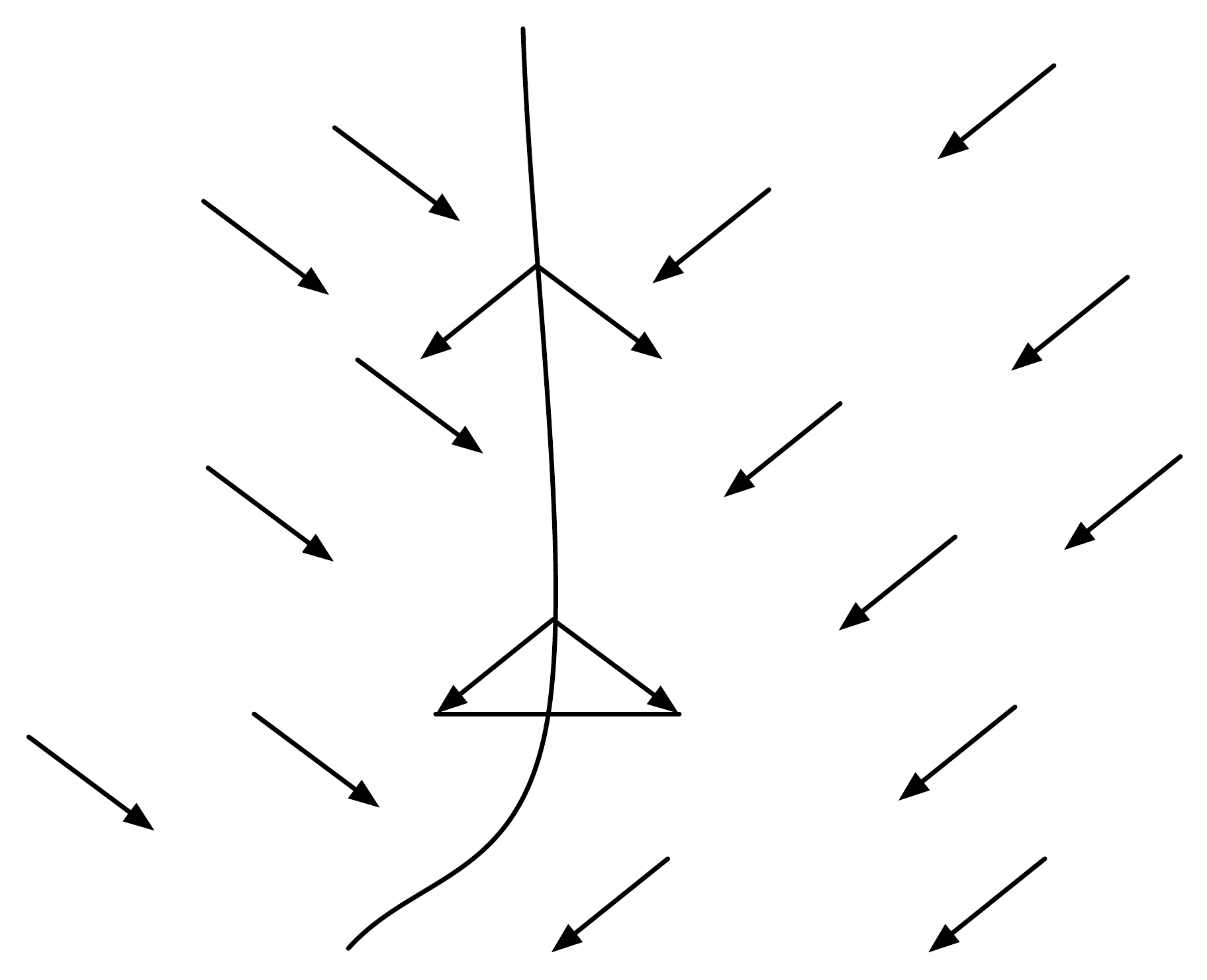}
  \caption{Intuitive meaning of Filippov Solutions on boundaries of partitions.}
\label{fig:control}
\end{figure}

It is also worth noting that $\mathcal{F}$ cannot be lower semi-continuous at any point $x \in X_i \cap X_j$, $ (i,j) \in \tilde{I}$, on a boundary, since $\mathcal{F}$ is not a one point set. For $T>0$, let $S_T$ denote either $[ 0,T)$ or $[ 0,T ]$. By a Carath\'eodory solution of differential inclusion \eqref{eq1} at $\zeta_0 \in \mathcal X$, we understand an absolutely continuous function $S_T \to \mathcal{X};t \mapsto \zeta(t)$ which solves the following Cauchy problem
\begin{equation} \label{eq7}
\dot{\zeta}(t) \in \mathcal{F}(\zeta(t)) \quad \text{a.e.,} \quad \zeta(0)=\zeta_0
\end{equation}
A Filippov solution to differential inclusion \eqref{eq1} at $\zeta_0 \in \mathcal X$ is a solution to \eqref{eq7} with $\mathcal{F}$ supplanted by $\mathcal{F}^c$ \cite{filip}. Intuitively, the concept of Filippov solutions implies that the velocity vector of a switched system exhibit a convex combination of velocity vectors in the neighborhood of a discontinuity or a boundary (see Fig.~\ref{fig:control}).

We recall the following facts from the theory of differential inclusions. Let $W$ be some non-negative function defined on $Graph(\mathcal{F})=\bigcup_{x \in \mathcal{X}} \{x\} \times \mathcal{F}(x) \subset \mathcal{X}^2$. We shall say that a function $\Phi : \mathcal{X} \to \mathbb{R}_{\ge 0}$ is a Lyapunov function for $\mathcal{F}$ with respect to $W$ if for all $x \in \mathcal{X}$ and some $v \in \mathcal{F}(x)$ the following ``Lyapunov property" holds
\begin{equation} \label{eq8}
D_+ \Phi (x)(v)+W(x,v) \le 0
\end{equation}
The next proposition asserts that under some mild conditions the switched system $\mathcal{S}$ allows for Filippov solutions.
\begin{prop}
Assume $F^c(\mathcal{X})$ is bounded. Then, at any point $x \in \mathcal{X}$, there exists a Filippov solution defined on ${S}_T$ to \eqref{eq1}. Furthermore, if it holds that
\begin{equation}\label{eq9}
\forall x \in \mathcal{X}, \quad \mathcal{F}^c(x) \cap T_\mathcal{X}(x) \ne \emptyset
\end{equation}
the solution exists on $S_\infty$.
\end{prop}
\begin{pf}
Because $\mathcal{F}^c$ is bounded, closed, convex, non-empty, and (from Proposition~\ref{p1}) upper semicontinuous everywhere on $\mathcal{X}$, from Corollary 1 p. 77 \cite{filip}, it can be concluded that there exists a Filippov solution to \eqref{eq1} on some bounded interval $S_T$ for any initial condition $x_0 \in \mathcal{X}$. Since $\mathcal{X}$ is compact, $F^c(\mathcal{X})$ is bounded, and \eqref{eq9} is satisfied, from Theorem 1 in p. 180 \cite{aubin} it follows that for all $x_0 \in \mathcal{X}$ there exists a solution to the differential inclusion $\dot x \in \mathcal{F}^c(x)$ with $x_0 \in \mathcal{X}$ defined on $\mathbb{R}_{\ge 0}$ which remains in $\mathcal{X}$.
\end{pf}

Now, we are ready to posit a stability condition for the set valued map $\mathcal{F}^c$, which is a consequence of applying Theorem 8.4 in p. 176, \cite{smir}.
\begin{prop} \label{p4}
Suppose $0 \in \mathcal{F}^c(0)$. If there exist $ \epsilon > 0$ and continuous positive definite functions $V:\mathbb{R}^n \to \mathbb{R}_{\ge 0}$ and $W:\mathbb{R}^n \to \mathbb{R}_{\ge 0}$ such that for each $x \in \mathcal{B}_\epsilon^n$
\begin{equation} \label{eq10}
D_+ V(x)(v)+W(x)\le 0 \quad \text{for all} \quad v \in \mathcal{F}^c(x),
\end{equation}
then the equilibrium point $0$ is asymptotically stable.

\subsection{Sum-of-Squares Programming}

\end{prop}
Recall that if there exists an SOS decomposition for $p(x)$, then it follows that $p(x)$ is non-negative. Unfortunately, the converse does not hold in general; that is, there exist non-negative polynomials which do not have an SOS decomposition. An epitome of this class of non-negative polynomials is the Motzkin's polynomial \cite{motz} given by
\begin{equation*}
p(x) = 1 - 3 x_1^2 x_2^2 + x_1^2 x_2^4 + x_1^4 x_2^2
\end{equation*}
which is non-negative for all $x \in \mathbb{R}^2$. This imposes, more or less, some sort of conservatism when utilizing SOS based methods. The next lemma gives an interesting formulation to the SOS decomposition problem.
\begin{lemm}[\cite{Choi95}] \label{p5}
A polynomial $p(x)$ of degree $2d$ belongs to $\mathcal{P}_{sos}$ if and only if there exist a positive semi-definite matrix $Q$ (known as the Gram matrix) and a vector of monomials $Z(x)$ which contains all monomial of $x$ of degree $\le d$ such that $p(x)=Z^T(x)QZ(x)$.
\end{lemm}
In \cite{chesi99}, Chesi \textit{et. al.} evinced that testing whether a
polynomial is SOS can be formulated as a set of LMI feasibility tests. Subsequently, Parrilo \cite{pari03} demonstrated that the answer to the query that whether a given polynomial $p(x)$ is SOS or not can be investigated via semi-definite programming methodologies.
\begin{lemm}[\cite{pari03}]
Given a finite set $\{ p_i \}_{i =0}^{m} \in  \mathcal{P}$, the existence of a set of scalars $\{ a_i \}_{i=1}^{m} \in \mathbb{R}$ such that
\begin{equation}
p_0 + \sum_{i=1}^m a_i p_i \in \mathcal{P}_{sos}
\end{equation}
is an LMI feasibility problem.
\end{lemm}
The subsequent lemma formalizes the problem of constrained positivity of polynomials which is a direct result of applying Stengle's Positivstellensatz method \cite{posi}.
\begin{lemm}[\cite{chesi}] \label{chesip}
Let $\{a_i\}_{i=1}^k$ and $\{b_i\}_{i=1}^l$ belong to $\mathcal{P}$, then
\begin{eqnarray}
p(x) \ge 0 \quad &\forall x \in \mathbb{R}^n: a_i(x)=0, \, \forall i=1,2,...,k & \nonumber \\
& \text{and} \quad b_j(x) \ge 0, \, \forall j=1,2,...,l&
\end{eqnarray}
is satisfied, if the following holds
\begin{eqnarray} \label{chesieq}
&\exists r_1,r_2,\ldots,r_k \in \mathcal{P} \quad \text{and} \quad \exists s_0,s_1,\ldots,s_l \in \mathcal{P}_{sos}& \nonumber \\
&p = \sum_{i=1}^k r_i a_i +\sum_{i=1}^l s_i b_i +s_0&
\end{eqnarray}
\end{lemm}
\begin{lemm} [\cite{chesi}] \label{spos}
The multivariable polynomial $p(x)$ is strictly positive ($p(x)>0 \quad \forall x \in \mathbb{R}^n$), if there exists a $\lambda > 0$ such that
\begin{equation}
\big( p(x) - \lambda \big) \in \mathcal{P}_{sos}
\end{equation}
\end{lemm}
At this point, we are prepared to delineate the main contributions of this paper.

\section{Main results} \label{secmain2}

In this section, we consider the stability problem of a class of nonlinear switched system defined on compact sets with Filippov solutions. Subsequently, we present a theorem for robust asymptotic stability of switched systems in the presence of polytopic uncertainty. Finally, we bring forward sufficient conditions for stability using SOS techniques.

\subsection{Asymptotic Stability Conditions for Switched Systems} \label{secmain}

Consider the switched system $\mathcal{S}$ and let \eqref{eq2} describe the Filippov solutions of $\mathcal{S}$. It is assumed that $0$ is an interior point of $\mathcal{X}$, and that it is located on some boundary of partitions. Note that $0 \in \mathcal{F}^c(0)$; hence, $0$ is an equilibrium. Define $\{ V_i(x) \}_{i \in I}$ a family of positive definite and continuously differentiable ($\mathcal{C}^1$) functions ($V_i : U_i \to \mathbb{R}_{\ge 0}$). We also define a set valued map $\Psi (x)$ associated with $\{ V_i(x) \}_{i \in I}$ as
\begin{equation} \label{phidef}
\Psi : \mathcal{X} \to 2^{\mathbb{R}}; x \mapsto \left\lbrace z \in \mathbb{R} \mid z=V_i(x) \quad \text{if} \quad x \in X_i \right\rbrace.
\end{equation}
We refer to $\Psi$ as a switched Lyapunov function. In general, this function cannot be continuously differentiable, since it is not a singleton for all $ x \in bd(X_i) \cup bd(X_{i+1}) \cup \cdots \cup bd(X_j)$.

\begin{prop} \label{cont}
If $V_i(x)=V_j(x)$ for all $x \in X_i \cap X_j$ and all $(i,j) \in \tilde{I}$, then $\Psi(\cdot)$ is real, single-valued ($\Psi : \mathcal{X} \to \mathbb{R}$), and locally Lipschitzean.
\end{prop}
Notice that, Proposition~\ref{cont} does not impose any constraint on the structure of $\{ V_i(x) \}_{i \in I}$, e.g. homogenous or quadratic forms as was done in \cite{lewis}. This considerably mitigates the conservatism in finding the family of Lyapunov functions $\{ V_i(x) \}_{i \in I}$. Once this family of functions is (somehow) found, one can directly construct the switched Lyapunov function.

\begin{prop} \label{acont}
Suppose
\begin{itemize}
\item[{I)}] $ \left \langle \frac{\partial{V_i(x)}}{\partial{x}} , F_i(x) \right \rangle < 0$ for all $x \in X_i \setminus \{0\}$ and all $i \in I$,
\item[{II)}] $\left \langle \frac{\partial{V_i(x)}}{\partial{x}} , F_j(x) \right \rangle < 0$ for all $x \in X_i \cap X_j \setminus \{0\}$ and all $(i,j) \in \tilde{I}$.
\end{itemize}
Then there exists a continuous positive definite function $W:\mathbb{R}^n \to \mathbb{R}_{\ge 0}$ such that
\begin{itemize}
\item[{III)}] $\left \langle \frac{\partial{V_i(x)}}{\partial{x}} , F_i(x) \right \rangle < -W(x)$ for all $x \in X_i \setminus \{0\}$ and \mbox{all $i \in I$},
\item[{IV)}] $\left \langle \frac{\partial{V_i(x)}}{\partial{x}} , F_j(x) \right \rangle < -W(x)$ for all $x \in X_i \cap X_j \setminus \{0\}$ and all $(i,j) \in \tilde{I}$.
\end{itemize}
\end{prop}
\begin{pf}
Suppose for each $X_i$, with $i \in I$, there exists an open neighborhood $T_i$ of $X_i$ such that condition (I) holds due to the compactness of $X_i$. Then, the collection of such open neighborhoods $\{ T_i \}_{i \in I}$ is an open cover of $\mathcal{X}$ such that $\mathcal{X} \subseteq \bigcup_{i \in I} T_{i}$. Therefore, there exists a partition of unity subordinate to the cover $\{ T_i \}_{i \in I}$; i.e, a family of continuous functions $\big \{ \psi_{i} : \mathcal{X} \to [0,1] \big \}_{i \in I}$ with $\mathnormal{supp}(\psi_{i}) \subset T_i$ such that for any point $x \in G$, there is a neighborhood of $x$ where all but finite number of functions $\{\psi_{i}\}_{i \in I}$ are equal to $0$, and such that $\sum_{i \in I} \psi_{i}(x) =1$. Thus, let $W_1(x)=- \sum_{i \in I} \psi_i(x) \langle \frac{\partial{V_i(x)}}{\partial{x}} , F_i(x) \rangle$ which satisfies (III). \\
In a similar manner, for all $Y_l = X_i \cap X_j$ with $(i,j) \in \tilde{I}$ and $l \in L=\{1,2,...,M\}$ (where $M$ is the number of members in $\tilde{I}$), there exist open neighborhoods $\mathcal{Y}_l$ whose collection ($\{ \mathcal{Y}_l \}_{l \in L}$) is an open cover to the closed set $\mathcal{X}' \subset \mathcal{X}$. Because $\mathcal{X}'$ is a closed subset of $\mathcal{X}$, $\mathcal{X}'$ is also compact. So, there exists a partition of unity subordinate to the cover $\{ \mathcal{Y}_l \}_{l \in L}$ characterized by $\big \{ \phi_{l} : \mathcal{Y}_l \to [0,1] \big \}_{l \in L}$. At this point, it suffices to let $W_2(x)=- \sum_{l \in L} \phi_l(x) \Gamma_l(x)$ where $\Gamma_l(x)=\langle \frac{\partial{V_i(x)}}{\partial{x}} , F_j(x) \rangle$ if $x \in Y_l$, $l \in L$. Obviously, $W_2(x)$ satisfies (IV). Finally, we can select the map $W(x)=\max \{ W_1(x), W_2(x) \}$, and this completes the proof.
\end{pf}

The next proposition provides a Lyapunov-like stability theorem for the class of switched systems under study.
\begin{them} \label{mainp}
Let $\{ V_i(x) \}_{i \in I}$ be a family of $\mathcal{C}^1$ Lyapunov functions. The switched system $\mathcal{S}$ is asymptotically stable at the origin if the following conditions hold
\begin{equation} \label{vpos}
V_i(x)>0 \qquad \forall x \in X_i \setminus \{0\}
\end{equation}
\begin{equation} \label{s11}
\left \langle \frac{\partial{V_i(x)}}{\partial{x}} , F_i(x) \right \rangle < 0 \qquad \forall x \in X_i \setminus \{0\}
\end{equation}
for all $i \in I$,
\begin{equation} \label{s2}
\left \langle \frac{\partial{V_i(x)}}{\partial{x}} , F_j(x) \right \rangle < 0 \qquad \forall x \in X_i \cap X_j \setminus \{0\}
\end{equation}
\begin{equation} \label{vcon}
V_i(x)=V_j(x) \qquad \forall x \in X_i \cap X_j \setminus \{0\}
\end{equation}
for all $(i,j) \in \tilde{I}$.
\end{them}
\begin{pf}
The proof follows the same lines as that of Proposition 10 in \cite{lewis}. It is necessary to show that conditions I-IV in Proposition~\ref{acont} holds true. From \eqref{phidef},\eqref{vpos},\eqref{vcon}, and Proposition~\ref{cont}, we conclude that there exists a continuous, locally Lipschitzean, single-valued, and positive definite function $\Psi(x)$. Subsequently, from \eqref{s11},\eqref{s2} and Proposition~\ref{acont} it follows that there exists a positive definite function $W(x)$ satisfying III and IV.\\
Given a set of $\mathcal{C}^1$ functions $\{ V_i (x) \}_{i \in I}$ and from the definition of construction covering, it follows that for any $x \in G$ and $v \in \mathbb{R}^n$, there is $a > 0$ such that $x + a t \in X_j$ for any $t \in [0, a)$. On the other hand, $D_+ V_j (x)(v) = \liminf _{h \to 0_+} \frac{V_j(x+hv) - V_j(x)}{h} = \langle \frac{\partial V_j}{\partial x}(x), v \rangle$. Then, from Proposition~\ref{acont} and \eqref{s11} it follows that
$$
D_{+} V_i (x) \big( F_i(x) \big) + W(x) \le 0
$$
Consequently, for any $u \in \mathcal{F}^c (x)$ and real $\alpha_k$ such that $\sum_{k \in I} \alpha_k =1$ , we arrive at the following justification
\begin{eqnarray}
D_{+} \Psi (x) (u) &=& D_{+} V_i(x) (u) \nonumber \\
&=& D_{+} V_i(x) \big( \sum_{k \in I} \alpha_k F_k (x) \big) \nonumber \\
&=& \bigg \langle \frac{\partial{V_i}}{\partial{x}} , \big( \sum_{k \in I} \alpha_k F_k (x) \big) \bigg \rangle \nonumber \\
&=& \sum_{k \in I} \alpha_k \big \langle \frac{\partial{V_i}}{\partial{x}} , F_k (x) \big \rangle \nonumber \\
& \le & - \sum_{k \in I} \alpha_k W(x) \nonumber \\
& \le & -W(x) \nonumber
\end{eqnarray}
in which, we applied \eqref{s2}, condition IV and Proposition~\ref{acont}. Thus, $0$ is an asymptotically stable equilibrium.
\end{pf}

\subsection{Robust Asymptotic Stability of Switched Systems with Simplical Uncertainty}
\begin{figure}
\centering
  \includegraphics[width=0.4\textwidth]{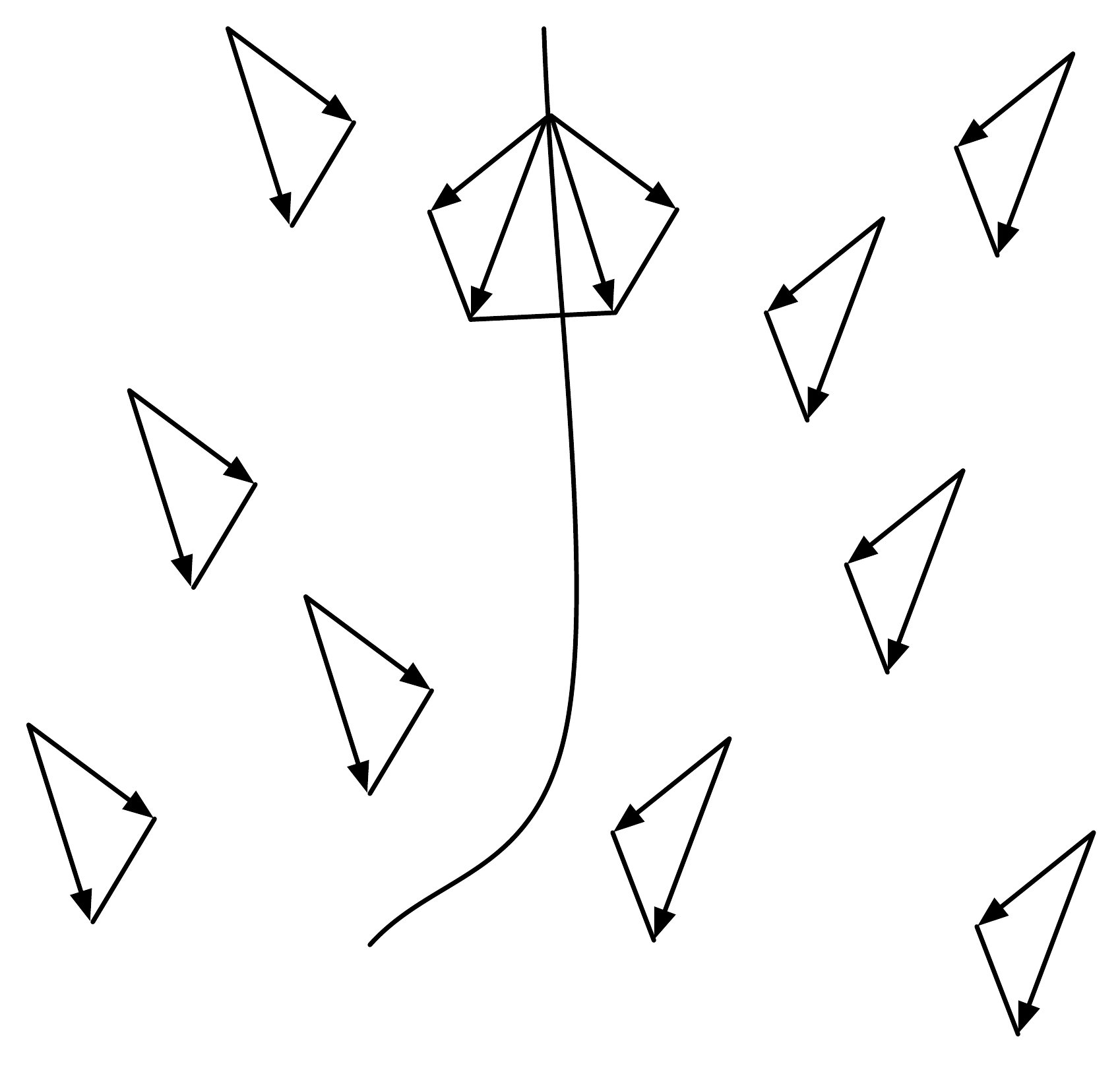}
  \caption{The dynamics of a switched system with simplical uncertainty in the neighborhood of a boundary.}
\label{fig:control2}
\end{figure}

At this stage, we extend our results to a class of switched systems with simplical uncertainty $\tilde{\mathcal{S}}=\lbrace \mathcal{X},\mathscr{X},I,\tilde{\mathscr{F}} \rbrace$ with $\tilde{\mathscr{F}}=\{ F_i(x,\theta^i) \}_{i \in I}$ and
\begin{equation} \label{defsu}
F_i(x,\theta^i) = \sum_{l=1}^{L_i} \theta_{il} f_{il}(x)
\end{equation}
where, $f_{il} : U_i \to \mathbb{R}^n, l=1,2,...,L_i$ ($U_i$ is an open neighborhood of $X_i$) are a family of smooth vector fields, and $\theta^i=\left (\theta_{i1},\theta_{i2},...,\theta_{iL_i} \right),\, i\in I$ are uncertain constant parameter vectors satisfying
\begin{eqnarray} \label{uncer}
\theta^i \in \Theta^i  \triangleq  \bigg \lbrace && \theta^i \in \mathbb{R}^{L_i} \mid \theta_{il} \ge 0,\, l=1,2,..,L_i, \, \nonumber \\
&& \text{and} \, \sum_{l=1}^{L_i} \theta_{il} = 1 \bigg \rbrace,
\end{eqnarray}
in which  each $\Theta^i$ is a simplex in $\mathbb R^{L_i}$. The presence of simplical uncertainty in the dynamics of a switched system can contribute to substantially discrepant motions than the ones dictated by each of the vector fields (see Fig. \ref{fig:control2}). Therefore, robust asymptotic stability of a switched system subject to uncertainty in the context of Filippov solutions seems to be a non-trivial problem. Fortunately,  with the results discussed in Section~\ref{secmain}, the following corollary regarding asymptotic stability for uncertain switched systems with Filippov solutions can be characterized.

\begin{coro} \label{pus}
Consider the switched system subject to simplical uncertainty $\tilde{\mathcal{S}}$. If there exists a family of Lyapunov functions $\{ V_i(x) \}_{i \in I}$ satisfying
\begin{equation} \label{l1}
V_i(x)>0 \qquad \forall x \in X_i \setminus \{0\}
\end{equation}
\begin{equation} \label{l2}
\big\langle \frac{\partial{V_i(x)}}{\partial{x}} , f_{il}(x) \big\rangle < 0 \qquad \forall x \in X_i \setminus \{0\}
\end{equation}
for all $i \in I$ and $l=1,2,...,L_i$,
\begin{equation} \label{l3}
\big\langle \frac{\partial{V_i(x)}}{\partial{x}} , f_{jl}(x) \big\rangle < 0 \qquad \forall x \in X_i \cap X_j \setminus \{0\}
\end{equation}
\begin{equation} \label{l4}
V_i(x)=V_j(x) \qquad \forall x \in X_i \cap X_j \setminus \{0\}
\end{equation}
for all $(i,j) \in \tilde{I}$, and $l=1,2,...,L_j$. Then, for all $\theta^i \in \Theta^i$ with $i \in I$, all Filippov solutions of $\tilde{\mathcal{S}}$ converge to origin asymptotically.
\end{coro}
\begin{pf}
This is a direct result of using Theorem~\ref{mainp}. Conditions \eqref{l1} and \eqref{l4} correspond to \eqref{vpos} and \eqref{vcon}, respectively. If \eqref{l2} holds for all $i \in I$ and $l=1,2,...,L_i$, then it follows that for all sets of unknown parameters $\theta^i$ satisfying \eqref{uncer}
\begin{eqnarray}
\sum_{l=1}^{L_i} \theta_{il} \big\langle \frac{\partial{V_i(x)}}{\partial{x}} , f_{il}(x) \big\rangle &=&
 \big\langle \frac{\partial{V_i(x)}}{\partial{x}} , \sum_{l=1}^{L_i} \theta_{il} f_{il}(x) \big\rangle \nonumber \\
 &=& \big\langle \frac{\partial{V_i(x)}}{\partial{x}} , F_i(x,\theta^i) \big\rangle < 0
\end{eqnarray}
which proves that \eqref{s11} holds. It can be analogously shown that if \eqref{l3} holds for all $(i,j) \in \tilde{I}$, and $l=1,2,...,L_j$, then it follows that \eqref{s2} is satisfied for all $(i,j) \in \tilde{I}$. Consequently, by Theorem~\ref{mainp}, all Filippov solutions of $\tilde{\mathcal{S}}$ converge to origin asymptotically.
\end{pf}

\subsection{Sufficient Conditions Based on SOS Programming}

Theorem~\ref{mainp} and Corollary~\ref{pus} present stability conditions for general nonlinear switched systems defined on compact sets. However, it is not clear how the family of functions $\{V_i(x)\}_{i \in I}$ is to be found. In order to supply Theorem~\ref{mainp} and Corollary \ref{pus} with computationally doable algorithms for constructing $\{V_i(x)\}_{i \in I}$, we assume that the vector fields associated with the switched systems are vectors of polynomials in the variable $x$, and that the switched systems are defined on semi-algebraic sets. Taking into account this assumption, we need computational efficient methods to check the positivity of a given polynomial over an specific set. The positivity test can be performed using two main approaches; i.e., the moments approach~\cite{doi:10.1137/S1052623400366802} and the SOS approach. In this study,  we use the SOS approach for which well-developed computational tools are available e.g. SOSTOOLS \cite{sos}.
Henceforth, we posit that each partition of $\mathcal{X}$ is described by a semi-algebraic set. For all $i \in I$, we have
\begin{equation} \label{epartition}
X_i = \left \{ x \in \mathcal{X} \mid \chi_i(x)=0, \,\, \xi_{ik}(x) \ge 0 \,\, \text{{for}} \,\, k \in N_i \right \},
\end{equation}
in which $\chi_i(x)$ and $\xi_{ik}$ with $i\in I$ and $k \in N_i$ belong to $\mathcal{P}(x)$. It can be readily deduced that $\chi_i(x)$ could take the form $\chi_i(x)=p_1^2(x)+p_2^2(x)+\cdots+p_n^2(x)$ for some $n>0$; hence, $\chi(x)=0$ implies that $p_1(x)=p_2(x)=\cdots=p_n(x)=0$. The boundary of partitions (switching surface) is a variety
\begin{equation} \label{eboundary}
X_i \cap X_j = \left \{ x \in \mathcal{X} \mid \chi_{ij}(x)=0 \right \}, \quad (i,j) \in \tilde{I}
\end{equation}
where $\chi_{ij}(x) \in \mathcal{P}(x)$ for all $(i,j) \in \tilde{I}$. The next theorem provides a set of SOS feasibility tests to construct $\{ V_i(x) \}_{i \in I}$, given the switched system $\mathcal{S}$ is asymptotically stable.
\begin{them} \label{mainpsos}
Let $\mathcal{S}$ be a polynomial switched system defined on semi-algebraic sets as described above. If there exist a family of polynomials $\{ V_i(x) \}_{i\in I}$ with $V_i(0)=0$ if $0 \in X_i$, $\phi_i(x) \in \mathcal{P}_{sos}(x)$, $\rho_i(x) \in \mathcal{P}(x)$, $\varepsilon_i(x) \in \mathcal{P}(x)$, $\{q_{ik}(x)\}_{k \in N_i} \subset \mathcal{P}_{sos}(x)$, $\{w_{ik}(x)\}_{k \in N_i} \subset \mathcal{P}_{sos}(x)$ with $i \in I$, and $p_{ij}(x) \in \mathcal{P}(x)$, $r_{ij}(x) \in \mathcal{P}(x)$ with $(i,j) \in \tilde{I}$, and two sets of positive scalars $\{\mu_i\}_{i \in I}$, and $\{\nu_{ij}\}_{(i,j) \in \tilde{I}}$ such that
\begin{eqnarray} \label{s1}
V_i(x)-\phi_i(x) &-&\varepsilon_i(x) \chi_i(x) \nonumber \\
&-&\sum_{k \in N_i} q_{ik}(x) \xi_{ik}(x)  \in \mathcal{P}_{sos}(x)
\end{eqnarray}
\begin{multline} \label{eee123}
- \big \langle \frac{\partial{V_i(x)}}{\partial{x}} , F_i(x) \big \rangle - \rho_i(x)\chi_i(x)  \\
-\sum_{k \in N_i} w_{ik}(x) \xi_{ik}(x) - \mu_i \in \mathcal{P}_{sos}(x)  
\end{multline}
for all $i \in I$, and
\begin{eqnarray} \label{slidsos1}
- \big \langle \frac{\partial{V_i(x)}}{\partial{x}} , F_j(x) \big \rangle - r_{ij}(x)\chi_{ij}(x)-\nu_{ij} \in \mathcal{P}_{sos}(x)
\end{eqnarray}
\begin{equation} \label{conte}
V_i(x)+p_{ij}(x)\chi_{ij}(x) =V_j(x)
\end{equation}
for all $(i,j) \in \tilde{I}$. Then, the equilibrium $0$ is asymptotically stable.
\end{them}
\begin{pf}
We need to apply Theorem~\ref{mainp}. Since $\mathcal{S}$ is defined on semi-algebraic sets the partitions and boundaries are given by \eqref{epartition} and \eqref{eboundary}, respectively. \eqref{conte} assures that \eqref{vcon} is satisfied; thus, Proposition~\ref{cont} holds. Since each boundary is a variety, condition \eqref{s2} can be reformulated using Lemma~\ref{chesip}  (see \eqref{chesieq} wherein $b_i=0$, $\rho_i(x)=r(x) \in \mathcal{P}(x)$, $p = - \big \langle \frac{\partial{V_i(x)}}{\partial{x}} , F_i(x) \big \rangle$, and noting that $s_0 \in \mathcal{P}_{sos}(x)$) and Lemma~\ref{spos}; consequently,  \eqref{slidsos1} is attained. Demonstrating that \eqref{eee123} corresponds to \eqref{s11} can be done in a similar fashion.  Furthermore, inasmuch as each partition is defined by a semi-algebraic set, \eqref{vpos} is analogous to
$$
V_i(x)>0, \,\, \text{if} \,\, \chi_i(x)=0 \,\, \text{and} \,\, \xi_{ik}(x) \ge 0 \,\, \text{{for}} \,\, k \in N_i
$$
Hence, using Lemma~\ref{chesip} and the generalized S-procedure \cite{slemma}, one can obtain \eqref{s1}. The term $-\phi_i(x)$ is added to ensure that each $V_i(x)$ is positive definite for all $x \in X_i$.
\end{pf}

It is worth noting that condition \eqref{slidsos1} can be further relaxed by just considering those boundaries possessing attractive Filippov solutions (instead of checking \eqref{slidsos1} for all $(i,j) \in \tilde{I}$). One can infer the existence of an attractive Filippov solution by checking
\begin{equation}
\big \langle \frac{\partial{\gamma_{ij}(x)}}{\partial{x}} , F_i(x) \big \rangle \big \langle \frac{\partial{\gamma_{ij}(x)}}{\partial{x}} , F_j(x) \big \rangle < 0 \quad \forall x \in X_i \cap X_j
\end{equation}
or, in terms of an SOS decomposition problem, if the following holds
\begin{multline} \label{attract}
-\big \langle \frac{\partial{\gamma_{ij}(x)}}{\partial{x}} , F_i(x) \big \rangle  \big \langle \frac{\partial{\gamma_{ij}(x)}}{\partial{x}} , F_j(x) \big \rangle  \\
 -l(x)\gamma_{ij}(x)-\kappa_{ij} \in \mathcal{P}_{sos}(x)
\end{multline}
for some $l(x) \in \mathcal{P}(x)$ and some positive scalar $\kappa_{ij}$. In fact, if for some $(i_*,j_*) \in \tilde{I}$ the SOS problem \eqref{attract} is not feasible; then, one can refrain from checking \eqref{slidsos1} for index pair $(i_*,j_*) \in \tilde{I}$.

It should be noted that Theorem~\ref{mainpsos} only provides sufficient conditions. Indeed, given a polynomial switched system $\mathcal{S}$ defined on semi-algebraic sets, one can search for the corresponding candidate Lyapunov functions via semi-definite programming schemes; if the problem is feasible, then the switched system $\mathcal{S}$ is asymptotically stable.

Based on similar arguments for Theorem~\ref{mainpsos}, we can characterize an SOS representation for conditions in Corollary~\ref{pus}.

\begin{coro} \label{p6}
Let $\tilde{\mathcal{S}}$ be a switched system with simplical uncertainty and defined on semi-algebraic sets. If there exist a family of polynomials $\{ V_i(x) \}_{i\in I}$ with $V_i(0)=0$ if $0 \in X_i$, $\phi_i(x) \in \mathcal{P}_{sos}(x)$, $\rho_i(x) \in \mathcal{P}(x)$, $\varepsilon_i(x) \in \mathcal{P}(x)$, $\{q_{ik}(x)\}_{k \in N_i} \subset \mathcal{P}_{sos}(x)$, $\{w_{ik}(x)\}_{k \in N_i} \subset \mathcal{P}_{sos}(x)$ with $i \in I$, and $p_{ij}(x) \in \mathcal{P}(x)$, $r_{ij}(x) \in \mathcal{P}(x)$ with $(i,j) \in \tilde{I}$, and two sets of positive scalars $\{\mu_{il}\}_{i \in I}$ $(l=1,2,...,L_i)$, and $\{\nu_{ijl}\}_{(i,j) \in \tilde{I}}$ $(l=1,2,...,L_j)$  such that \eqref{s1} and
\begin{multline}
- \big \langle \frac{\partial{V_i(x)}}{\partial{x}} , f_{il}(x) \big \rangle -\rho_i(x)\chi_i(x)  \\ -\sum_{k \in N_i} w_{ik}(x) \xi_{ik}(x) - \mu_{il} \in \mathcal{P}_{sos}(x)  \\
\end{multline}
holds for all for all $i \in I$ and $l=1,2,...,L_i$,
\begin{equation} \label{slidsos}
- \big \langle \frac{\partial{V_i(x)}}{\partial{x}} , f_{jl}(x) \big \rangle -r_{ij}(x)\chi_{ij}(x)-\nu_{ijl} \in \mathcal{P}_{sos}(x)
\end{equation}
and \eqref{conte} holds for all for all $(i,j) \in \tilde{I}$ and $l=1,2,...,L_j$. Then, the origin is robustly asymptotically stable.
\end{coro}

\section{Simulation Analysis} \label{secexam}

\begin{figure*}[t]
\begin{center}
\resizebox{7cm}{!}{
  \includegraphics
    {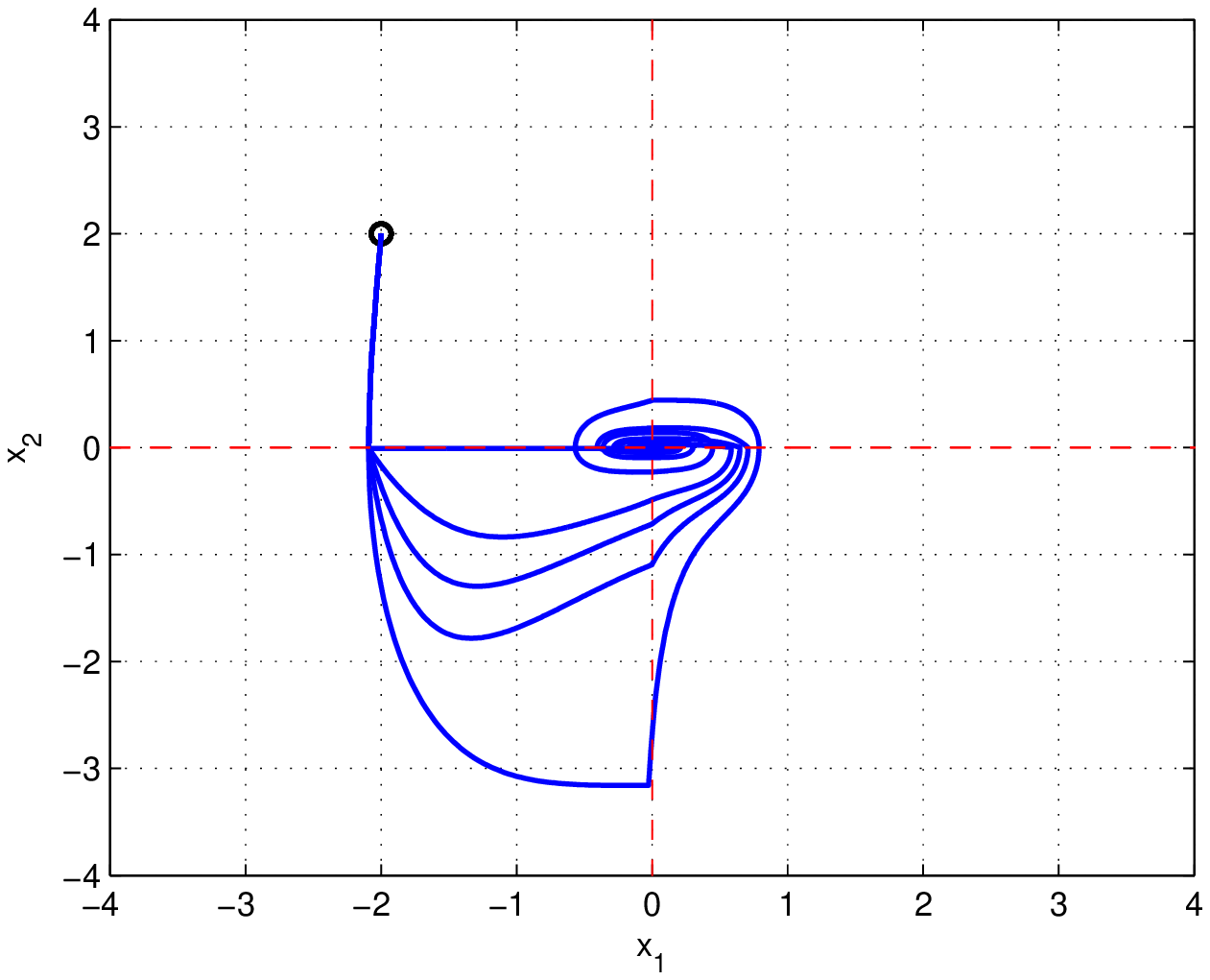}
}
\resizebox{7cm}{!}{
  \includegraphics
        {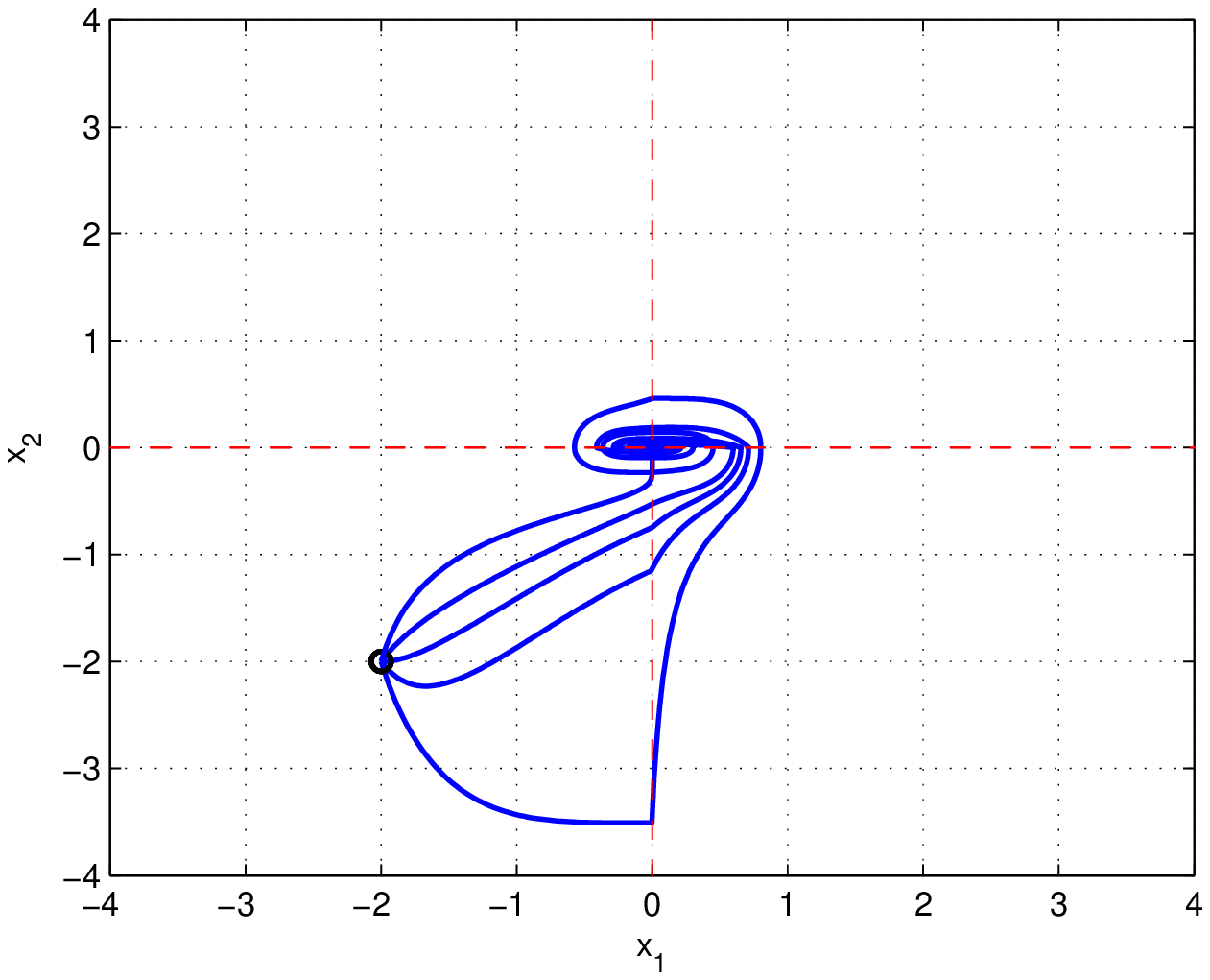}
}
\\
\resizebox{7cm}{!}{
  \includegraphics
        {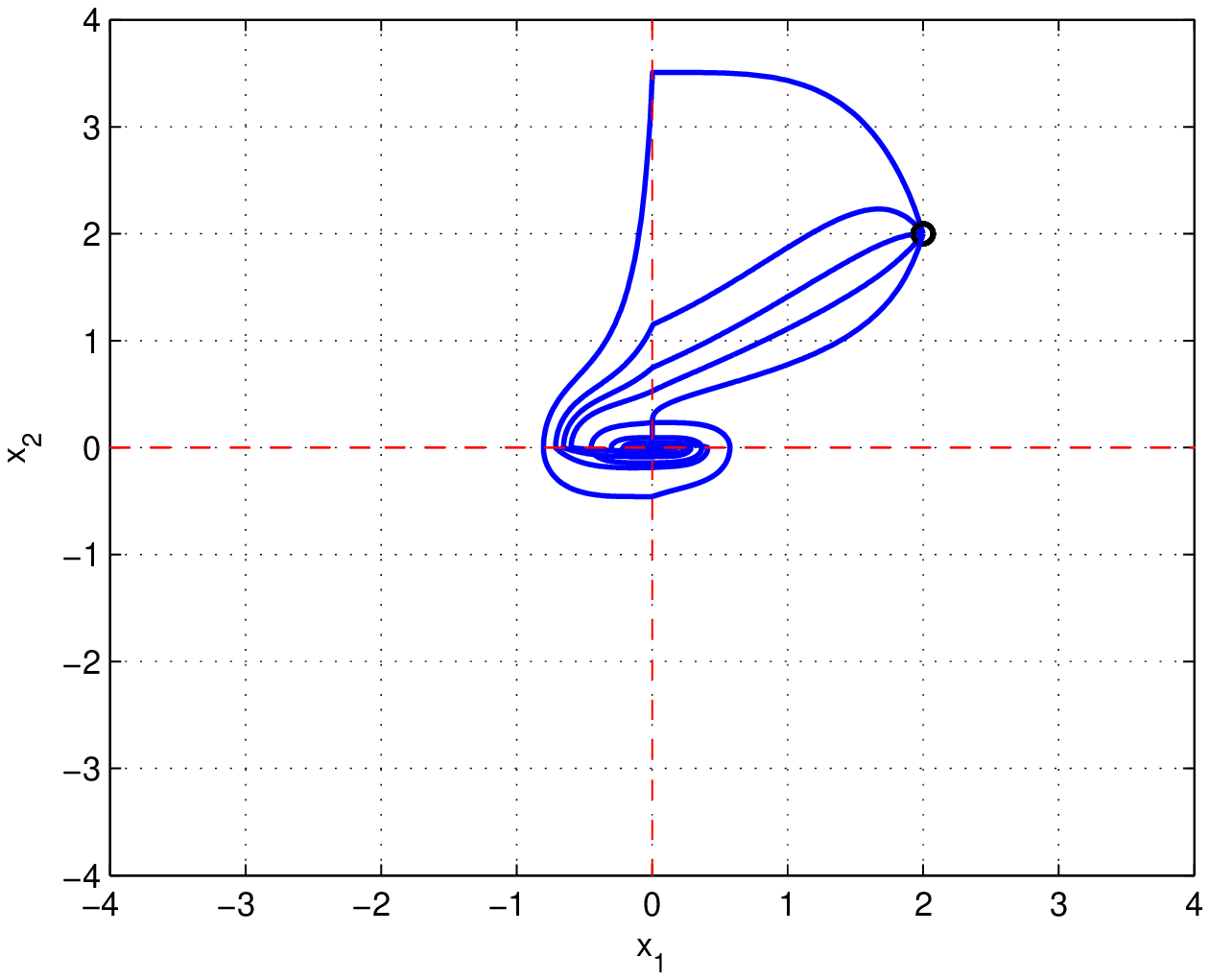}
}
\resizebox{7cm}{!}{
  \includegraphics
        {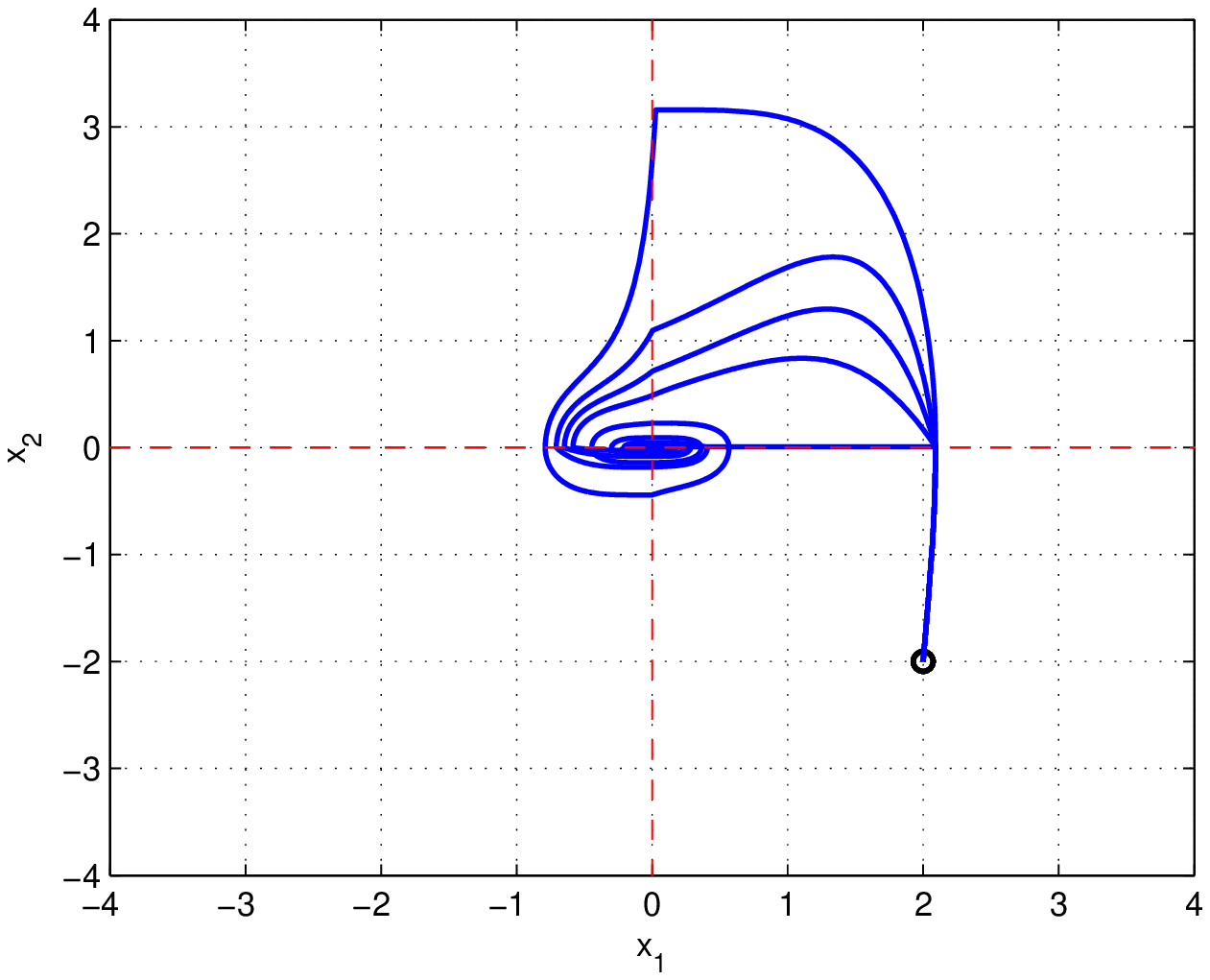}
}
  \caption{ Trajectories of the switched system with four starting points in different quadrants.}
  \label{figex2}
\end{center}
\end{figure*}

\begin{figure}[t]
\begin{center}

{
\resizebox{6cm}{!}{
  \includegraphics
        {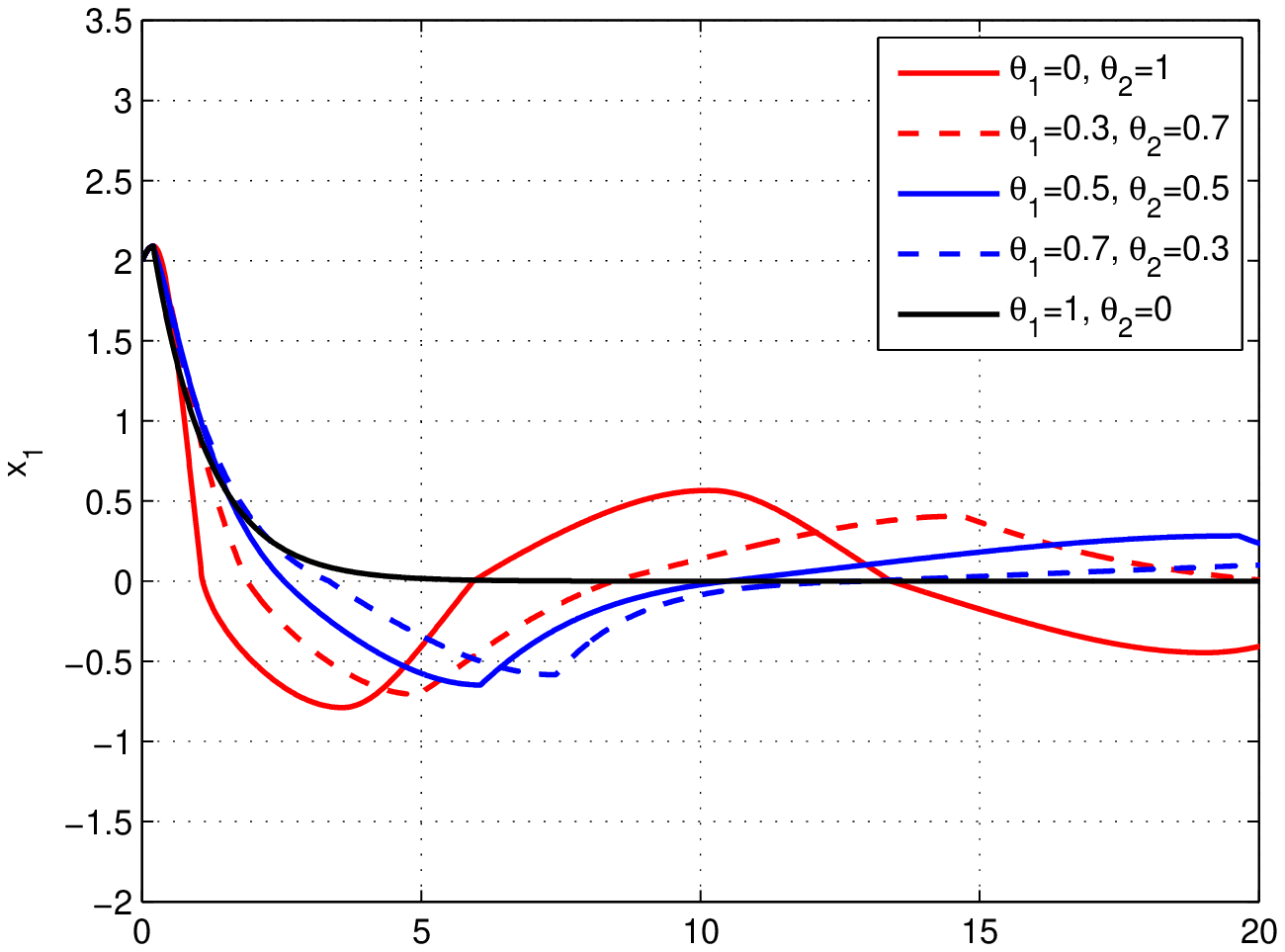}
}
\\
\resizebox{6cm}{!}{
  \includegraphics
        {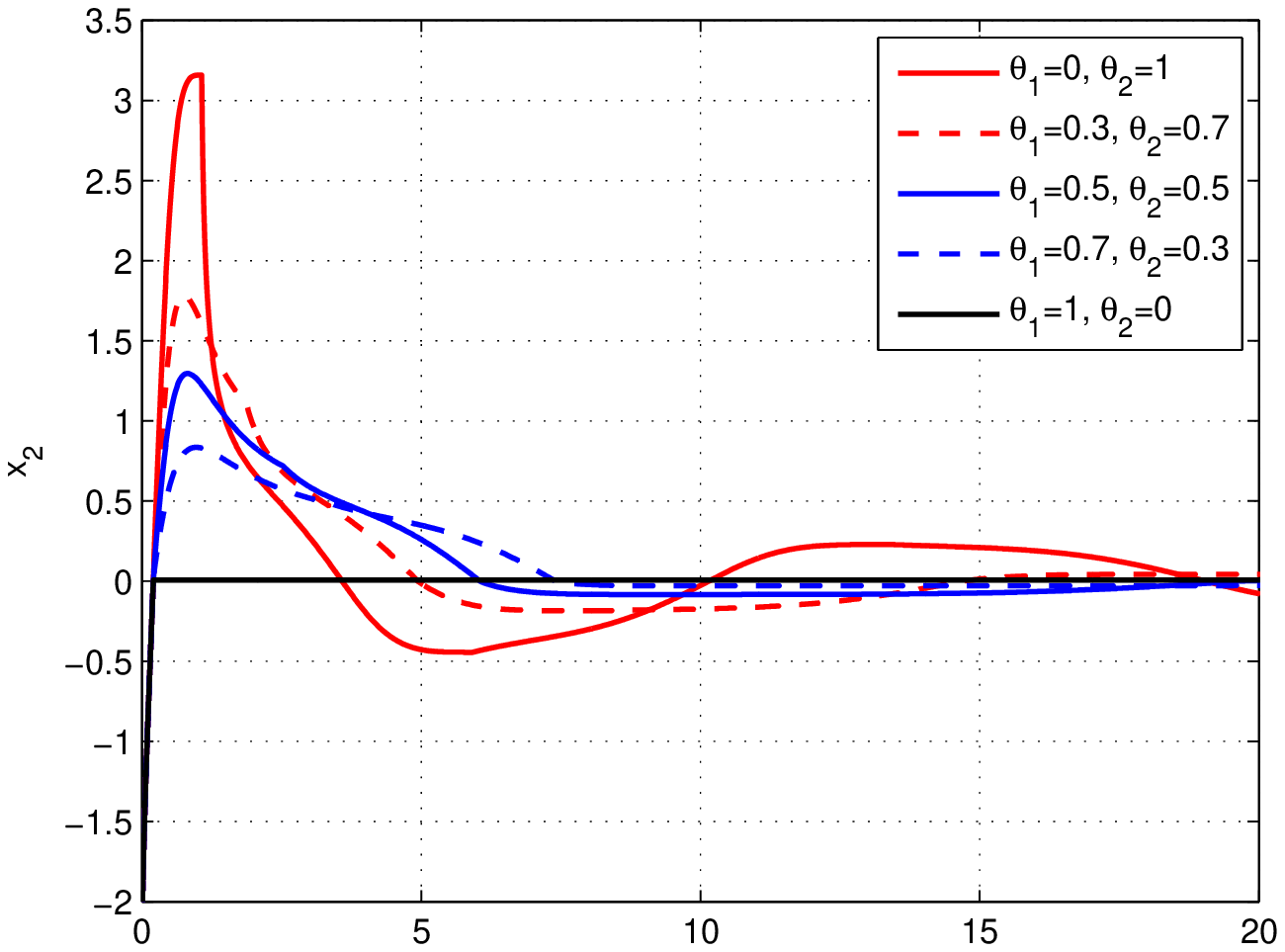}
}
\\
\resizebox{6cm}{!}{
  \includegraphics
        {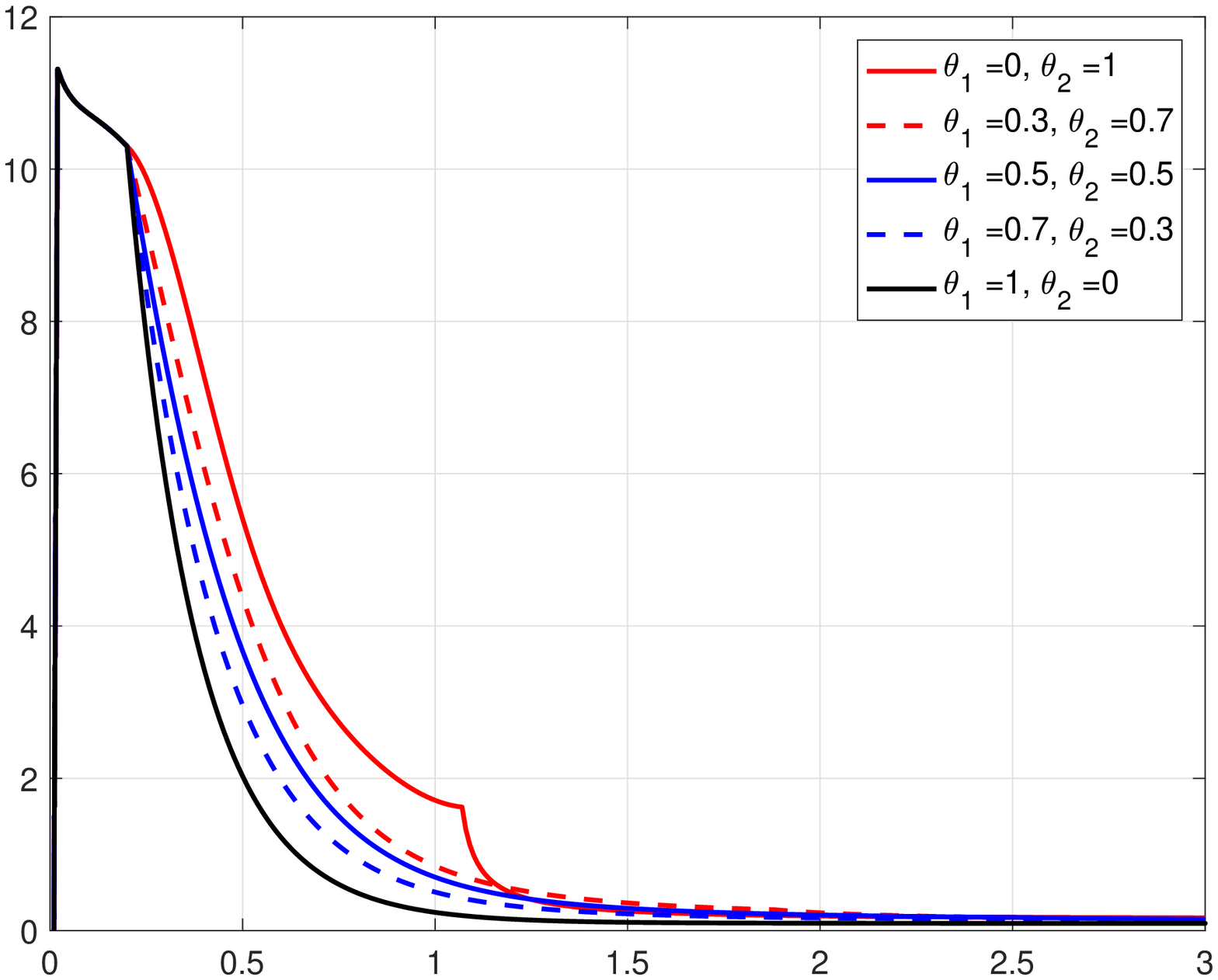}
}}

  \caption{ State responses for the switched system in the numerical example starting at $x_0=(2,-2)^T$ and the corresponding values of $\Psi(x)$.}
  \label{lasttt2}
\end{center}
\end{figure}

Consider an uncertain switched system where the dynamics is described by
\begin{eqnarray}
\dot x & \in & \tilde{\mathcal H} (x,\theta) \label{eq110} \\
\dot x & \in & co \left ( \tilde{\mathcal H} (x,\theta) \label{eq111} \right ),
\end{eqnarray}
where $H_i(x,\theta^i)=\sum_{l=1}^2 \theta_{il} h_{il}(x)$ with $\theta_{21}=0$ (obviously $\theta_{22}=1$) and $I=\{1,2\}$, the partitions defined as
\begin{eqnarray}
&X_1 = \lbrace x \in \mathbb{R}^2 \mid x_1 x_2 > 0 \rbrace&\\
&X_2 = \lbrace x \in \mathbb{R}^2 \mid x_1 x_2 < 0 \rbrace&,
\end{eqnarray}
and the local subsystems described by
\begin{eqnarray}
&h_{11}(x) =
\begin{bmatrix}
-x_1 \\
-x_2^3
\end{bmatrix}, \quad h_{12}(x) =
\begin{bmatrix}
-x_2 \\
x_1^3
\end{bmatrix}&\\
&h_{22}(x) =
\begin{bmatrix}
-0.5 x_2 \\
x_1^3 - x_2^3
\end{bmatrix}.&
\end{eqnarray}
Notice that partitions $X_1$ and $X_2$ are not convex. As demonstrated in Fig. \ref{figex2}, the simulations attest that the uncertain system is asymptotically stable for all values of $\theta_{11}$, and $\theta_{12}$. The five trajectories in each subfigure corresponds to parameter values of $\theta_{11}=0,0.3,0.5,0.7,1$. 
The following Lyapunov functions of degree $6$ (removing the terms with coefficients smaller than $10^{-4}$) were determined using SOSTOOLS in 4.6956 seconds on a personal computer with  Intel(R) Core(TM) 2 Due CPU T7500 @ 2.20GHz and 4 GB of RAM
\begin{eqnarray}
V_1(x) &=&  0.98847x_1^6 + 1.0184x_1^4  - 0.70253x_1^3x_2  \nonumber \\
&&+ 0.035508x_1x_2^3 + 0.0037x_2^4   \nonumber \\
&& + 0.0021x_1x_2 + 1.113x_1^2x_2^2 + 1.4904x_2^2 \nonumber \\
&&+ 0.95437;  \label{lex1} \\
V_2(x) &=& V_1(x) +0.001x_1^2x_2^2. \label{lex2}
\end{eqnarray}
Consequently from Corollary~\ref{p6}, one can conclude that the uncertain switched system with Filippov solutions given by inclusions \eqref{eq110} and \eqref{eq111} is robustly asymptotically stable.
We can associate with \eqref{lex1} and \eqref{lex2} a switched Lyapunov function as given in \eqref{phidef}. Figs.  \ref{lasttt2} illustrates the Filippov solutions of \eqref{eq110} and the corresponding evolution of the switched Lyapunov function $\Psi(x)$. Notice that for the sake of convenience of display in the figures, $\theta_{11}$ and $\theta_{12}$ are represented respectively as $\theta_1$ and $\theta_2$. It can be discerned from the figure that due to simplical uncertainty, the solutions of the switched system may possess different behaviors in terms of stability. However, the switched Lyapunov function is still positive definite and dwindling as the Filippov solutions converge to origin.

%
%

\section{Conclusion} \label{seccon}

A Lyapunov-like stability theorem for nonlinear switched systems with partitioned state-space and state-dependent switching is brought forward. This result was used to formulate conditions for robust asymptotic stability of switched systems with simplical uncertainty. Since the analysis is based on the theory of differential inclusions, the proposed stability analysis scheme includes Filippov solutions or sliding modes. Furthermore, for the case of polynomial switched systems defined on semi-algebraic sets, we provide a computationally efficient method  based on SOS programming to implement the suggested stability theorems. The validity of the proposed methodologies was examined through simulation analysis.

\section*{References}

\bibliographystyle{model1-num-names}
\bibliography{autosam}           



\appendix
\end{document}